\documentclass[11pt]{article}
\usepackage{amssymb}
\usepackage{amsmath}
\usepackage[all]{xy}

\def\Mod{\mathop{\rm Mod}\nolimits}

\def\Ext{\mathop{\rm Ext}\nolimits}

\def\Hom{\mathop{\rm Hom}\nolimits}
\def\Im{\mathop{\rm Im}\nolimits}

\def\GProj{\mathop{\rm GProj}\nolimits}
\def\GInj{\mathop{\rm GInj}\nolimits}
\def\Proj{\mathop{\rm Proj}\nolimits}
\def\Inj{\mathop{\rm Inj}\nolimits}

\title{\Large \bf Duality of Preenvelopes and Pure Injective Modules
}
\author{Zhaoyong Huang \\
{\footnotesize \it Department of Mathematics, Nanjing University,
Nanjing 210093, Jiangsu Province, P. R. China}
\\{\footnotesize \it (email: huangzy@nju.edu.cn)}}
\date{}
\begin{document}
\baselineskip=18pt \maketitle

\begin{abstract}
Let $R$ be an arbitrary ring and $(-)^+=\Hom _{\mathbb{Z}}(-,
\mathbb{Q}/\mathbb{Z})$ where $\mathbb{Z}$ is the ring of integers
and $\mathbb{Q}$ is the ring of rational numbers, and let
$\mathcal{C}$ be a subcategory of left $R$-modules and $\mathcal{D}$
a subcategory of right $R$-modules such that $X^+\in \mathcal{D}$
for any $X\in \mathcal{C}$ and all modules in $\mathcal{C}$ are pure
injective. Then a homomorphism $f: A\to C$ of left $R$-modules with
$C\in \mathcal{C}$ is a $\mathcal{C}$-(pre)envelope of $A$ provided
$f^+: C^+\to A^+$ is a $\mathcal{D}$-(pre)cover of $A^+$. Some
applications of this result are given.
\end{abstract}

\vspace{0.5cm}

{\bf 1. Introduction}

\vspace{0.2cm}

Throughout this paper, all rings are associative with identity. For
a ring $R$, we use $\Mod R$ (resp. $\Mod R^{op}$) to denote the
category of left (resp. right) $R$-modules.

(Pre)envelopes and (pre)covers of modules were introduced by Enochs
in [E], and are fundamental and important in relative homological
algebra. Following Auslander and Smal$\phi$'s terminology in [AuS],
for a finitely generated module over an artinian algebra, a
(pre)envelope and a (pre)cover are called a (minimal) left
approximation and a (minimal) right approximation, respectively.
Notice that (pre)envelopes and (pre)covers of modules are dual
notions, so the dual properties between them are natural research
topics. It has been known that most of their properties are indeed
dual ([AuS, E, EH, EJ2, GT] and references therein).

We write $(-)^+=\Hom _{\mathbb{Z}}(-, \mathbb{Q}/\mathbb{Z})$, where
$\mathbb{Z}$ is the ring of integers and $\mathbb{Q}$ is the ring of
rational numbers. For a ring $R$ and a subcategory $\mathcal{X}$ of
$\Mod R$ (or $\Mod R^{op}$), we write $\mathcal{X}^+=\{X^+\ |\ X\in
\mathcal{X}\}$. Enochs and Huang proved in [EH] the following
result, which played a crucial role in [EH].

\vspace {0.2cm}

{\bf Theorem 1.1} ([EH, Corollary 3.2]) {\it Let $R$ be a ring, and
let $\mathcal{C}$ be a subcategory of $\Mod R$ and $\mathcal{D}$ a
subcategory of $\Mod R^{op}$ such that $\mathcal{C}^+\subseteq
\mathcal{D}$ and $\mathcal{D}^+\subseteq \mathcal{C}$. If $f: A\to
C$ is a $\mathcal{C}$-preenvelope of a module $A$ in $\Mod R$, then
$f^+: C^+\to A^+$ is a $\mathcal{D}$-precover of $A^+$ in $\Mod
R^{op}$.}

\vspace{0.2cm}

However, the converse of Theorem 1.1 does not hold true in general
(see [EH, Example 3.6]). So a natural question is: when does the
converse of Theorem 1.1 hold true? In this paper, we will give a
partial answer to this question and prove the following

\vspace{0.2cm}

{\bf Theorem 1.2.} {\it  Let $R$ be a ring, and let $\mathcal{C}$ be
a subcategory of $\Mod R$ and $\mathcal{D}$ a subcategory of $\Mod
R^{op}$ such that $\mathcal{C}^+\subseteq \mathcal{D}$ and all
modules in $\mathcal{C}$ are pure injective. Then a homomorphism $f:
A\to C$ in $\Mod R$ with $C\in \mathcal{C}$ is a
$\mathcal{C}$-preenvelope of $A$ provided $f^+: C^+\to A^+$ is a
$\mathcal{D}$-precover of $A^+$ in $\Mod R^{op}$.}

\vspace{0.2cm}

This paper is organized as follows.

In Section 2, we give some terminology and some preliminary results.

Let $R$ and $S$ be rings and let $_SU_R$ be a given $(S,
R)$-bimodule. For a subcategory $\mathcal{X}$ of $\Mod S$ (or $\Mod
R^{op}$), we write $\mathcal{X}^*=\{X^*\ |\ X\in \mathcal{X}\}$,
where $(-)^*=\Hom(-, {_SU_R})$. In Section 3, we first prove that if
$\mathcal{C}$ is a subcategory of $\Mod S$ and $\mathcal{D}$ is a
subcategory of $\Mod R^{op}$ such that $\mathcal{C}^*\subseteq
\mathcal{D}$ and the canonical evaluation homomorphism $X\to X^{**}$
is a split monomorphism for any $X\in \mathcal{C}$, then a
homomorphism $f: A\to C$ in $\Mod S$ being a
$\mathcal{C}$-preenvelope of $A$ implies that $f^*: C^*\to A^*$ is a
$\mathcal{D}$-precover of $A^*$ in $\Mod R^{op}$. As a consequence
of this result we get easily Theorem 1.2. Then as applications of
Theorem 1.2 we get the following results. For a ring $R$, a
monomorphism $f: A\rightarrowtail C$ in $\Mod R$ with $C$ pure
injective is a pure injective (pre)envelope of $A$ provided $f^+:
C^+\twoheadrightarrow A^+$ is a pure injective (or cotorsion)
(pre)cover of $A^+$ in $\Mod R^{op}$. For a left and right artinian
ring $R$, a homomorphism $f: A\to P$ in $\Mod R$ is a projective
preenvelope of $A$ if and only if $f^+: P^+\to A^+$ is an injective
precover of $A^+$ in $\Mod R^{op}$. In particular, we prove that for
a left and right coherent ring $R$, an absolutely pure left
$R$-module does not have a decomposition as a direct sum of
indecomposable absolutely pure submodules in general. It means that
a left (and right) coherent ring has no absolutely pure analogue of
[M, Theorem 2.5], which states that for a left noetherian ring $R$,
every injective left $R$-module has a decomposition as a direct sum
of indecomposable injective submodules.

\vspace{0.5cm}

{\bf  2. Preliminaries}

\vspace{0.2cm}

In this section, we give some terminology and some preliminary
results for later use.

\vspace{0.2cm}

{\bf Definition 2.1.} ([E]) Let $R$ be a ring and $\mathcal{C}$ a
subcategory of $\Mod R$. The homomorphism $f: C\to D$ in $\Mod R$
with $C\in \mathcal{C}$ is said to be a {\it $\mathcal{C}$-precover}
of $D$ if for any homomorphism $g: C' \to D$ in $\Mod R$ with $C'\in
\mathcal{C}$, there exists a homomorphism $h: C'\to C$ such that the
following diagram commutes:
$$\xymatrix{ & C' \ar[d]^{g} \ar@{-->}[ld]_{h}\\
C \ar[r]^{f} & D}$$ The homomorphism $f: C\to D$ is said to be {\it
right minimal} if an endomorphism $h: C\to C$ is an automorphism
whenever $f=fh$. A $\mathcal{C}$-precover $f: C\to D$ is called a
{\it $\mathcal{C}$-cover} if $f$ is right minimal. Dually, the
notions of a {\it $\mathcal{C}$-preenvelope}, a {\it left minimal
homomorphism} and {\it a $\mathcal{C}$-envelope} are defined.

\vspace{0.2cm}

Let $R$ be a ring. Recall that a short exact sequence $0 \to A \to
B\to C \to 0$ in $\Mod R$ is called {\it pure} if the functor $\Hom
_R(M,-)$ preserves its exactness for any finitely presented left
$R$-module $M$, and a module $E\in \Mod R$ is called {\it pure
injective} if the functor $\Hom _R(-,E)$ preserves the exactness of
a short pure exact sequence in $\Mod R$ (cf. [GT, K]). Recall from
[K] that a subcategory $\mathcal{C}$ of $\Mod R$ is {\it definable}
if it is closed under direct limits, direct products and pure
submodules in $\Mod R$.

\vspace{0.2cm}

{\bf Lemma 2.2.} {\it ([K, Corollary 2.7]) The following statements
are equivalent for a definable subcategory $\mathcal{C}$ of $\Mod
R$.

(1) Every module in $\mathcal{C}$ is pure injective.

(2) Every module in $\mathcal{C}$ is a direct sum of indecomposable
modules.}

\vspace{0.2cm}

{\bf Lemma 2.3.} {\it (1) ([F, Theorem 2.1]) For a ring $R$, a
module $M$ in $\Mod R$ is flat if and only if $M^+$ is injective in
$\Mod R^{op}$.

(2) ([CS, Theorem 4]) A ring $R$ is left (resp. right) artinian if
and only if a module $A$ in $\Mod R$ (resp. $\Mod R^{op}$) being
injective is equivalent to $A^+$ being projective in $\Mod R^{op}$
(resp. $\Mod R$).}

\vspace{0.2cm}

As a generalization of projective (resp. injective) modules, the
notion of Gorenstein projective (resp. injective) modules was
introduced by Enochs and Jenda in [EJ1] as follows.

\vspace{0.2cm}

{\bf Definition 2.4.} ([EJ1]) Let $R$ be a ring. A module $M$ in
$\Mod R$ is called {\it Gorenstein projective} if there exists an
exact sequence:
$$\mathbb{P}:\cdots \to P_1 \to P_0 \to P^0 \to P^1 \to \cdots$$
in $\Mod R$ with all terms projective, such that $M=\Im (P_0 \to
P^0)$ and the sequence $\Hom _R(\mathbb{P}, P)$ is exact for any
projective left $R$-module $P$. Dually, the notion of Gorenstein
injective modules is defined.

\vspace{0.5cm}

{\bf  3. The duality between preenvelopes and precovers}

\vspace {0.2cm}


Let $R$ and $S$ be rings and let $_SU_R$ be a given
$(S,R)$-bimodule. We write $(-)^*=\Hom(-, {_SU_R})$. For a
subcategory $\mathcal{X}$ of $\Mod S$ (or $\Mod R^{op}$), we write
$\mathcal{X}^*=\{X^*\ |\ X\in \mathcal{X}\}$. For any $X\in \Mod S$
(or $\Mod R^{op}$), $\sigma _X: X \to X^{**}$ defined by $\sigma
_X(x)(f)=f(x)$ for any $x\in X$ and $f\in X^*$ is the canonical
evaluation homomorphism.







The following lemma plays a crucial role in proving the main result.

\vspace{0.2cm}

{\bf Lemma 3.1.} {\it Let $\mathcal{C}$ be a subcategory of $\Mod S$
and $\mathcal{D}$ a subcategory of $\Mod R^{op}$ such that
$\mathcal{C}^*\subseteq \mathcal{D}$ and $\sigma _X$ is a split
monomorphism for any module $X\in \mathcal{C}$. For a homomorphism
$f: A\to C$ in $\Mod S$ with $C\in \mathcal{C}$, if $f^*: C^*\to
A^*$ is a $\mathcal{D}$-precover of $A^*$ in $\Mod R^{op}$, then $f:
A\to C$ is a $\mathcal{C}$-preenvelope of $A$.}

\vspace{0.2cm}

{\it Proof.} Let $f: A\to C$ be in $\Mod S$ with $C\in \mathcal{C}$
such that $f^*: C^*\to A^*$ is a $\mathcal{D}$-precover of $A^*$ in
$\Mod R^{op}$. Assume that $A\in \Mod S$, $X\in \mathcal{C}$ and
$g\in \Hom_S(A,X)$. Then $X^*\in \mathcal{C}^*\subseteq \mathcal{D}$
and there exists $h\in \Hom_{R^{op}}(X^*,C^*)$ such that the
following diagram commutes:
$$\xymatrix{C^* \ar[r]^{f^*} & A^*\\
& X^* \ar[u]_{g^*}\ar@{-->}[lu]^{h}}$$ Then $g^*=f^*h$ and
$g^{**}=h^*f^{**}$.

We have the following diagram with each square commutative:
$$\xymatrix{X \ar[r]^{\sigma_X} & X^{**}\\
A \ar[r]^{\sigma_A}\ar[u]^{g}\ar[d]_{f} & A^{**}\ar[u]_{g^{**}}\ar[d]^{f^{**}}\\
C \ar[r]^{\sigma_C} & C^{**}}$$ Then $g^{**}\sigma_A=\sigma_X g$ and
$f^{**}\sigma_A=\sigma_C g$. By assumption $\sigma _X$ is a split
monomorphism, so there exists $\alpha\in\Hom_S(X^{**},X)$ such that
$\alpha\sigma_X=1_X$, and hence we have that $g=1_X
g=(\alpha\sigma_X)g=\alpha g^{**}\sigma_A=\alpha
(h^*f^{**})\sigma_A=(\alpha h^*\sigma_C)f$, that is, we get a
homomorphism $\alpha h^*\sigma_C:C\to X$ in $\Mod S$ such that the
following diagram commutes:
$$\xymatrix{A \ar[r]^f \ar[d]_g & C\ar@{-->}[ld]^{\alpha h^*\sigma_C}\\
X. &}$$ Thus $f$ is a $\mathcal{C}$-preenvelope of $A$.
\hfill{$\square$}

\vspace{0.2cm}

From now on, $R$ is an arbitrary ring. For a subcategory
$\mathcal{X}$ of $\Mod R$ (or $\Mod R^{op}$), we write
$\mathcal{X}^+=\{X^+\ |\ X\in \mathcal{X}\}$. The main result in
this paper is the following







\vspace {0.2cm}

{\bf Theorem 3.2.} {\it Let $\mathcal{C}$ be a subcategory of $\Mod
R$ and $\mathcal{D}$ a subcategory of $\Mod R^{op}$ such that
$\mathcal{C}^+\subseteq \mathcal{D}$ and all modules in
$\mathcal{C}$ are pure injective. Then a homomorphism $f: A\to C$ in
$\Mod R$ with $C\in \mathcal{C}$ is a $\mathcal{C}$-(pre)envelope of
$A$ provided $f^+: C^+\to A^+$ is a $\mathcal{D}$-(pre)cover of
$A^+$ in $\Mod R^{op}$.}

\vspace{0.2cm}

{\it Proof.} First note that
$\Hom_{R^{op}}(-, R^+)\cong (-)^+$ by the adjoint isomorphism
theorem. Now let $X$ be a module in $\mathcal{C}$. Then $C$ is pure
injective by assumption, and so $\sigma_X: X\to X^{++}$ is a split
monomorphism by [GT, Theorem 1.2.19]. So the assertion follows from
Lemma 3.1 and [EH, Corollary 3.2(2)]. \hfill{$\square$}

\vspace{0.2cm}

In the rest of this section, we will give some applications of
Theorem 3.2.

By [K, Example 3.16], any module in $\Mod R$ has a pure injective
envelope. It is obvious that the pure injective (pre)envelope of any
module is monic. Recall from [EJ2] that a module $N\in \Mod R^{op}$
is called {\it cotorsion} if $\Ext_{R^{op}}^1(F,N)=0$ for any flat
right $R$-module $F$. For any $M\in \Mod R$, $M^+$ is pure injective
by [EJ2, Proposition 5.3.7], and hence cotorsion by [EJ2, Lemma
5.3.23]. Then by Theorem 3.2, we immediately have the following

\vspace{0.2cm}

{\bf Corollary 3.3.} {\it A monomorphism $f: A\rightarrowtail C$ in
$\Mod R$ with $C$ pure injective is a pure injective (pre)envelope
of $A$ provided $f^+: C^+\twoheadrightarrow A^+$ is a pure injective
(or cotorsion) (pre)cover of $A^+$ in $\Mod R^{op}$.}

\vspace{0.2cm}

Recall that $R$ is called {\it left pure semisimple} if every left
$R$-module is a direct sum of finitely generated modules, or
equivalently, every left $R$-module is pure injective.
As an immediate consequence of Theorem 3.2, we have the following

\vspace{0.2cm}

{\bf Corollary 3.4.} {\it Let $R$ be a left pure semisimple ring,
and let $\mathcal{C}$ be a subcategory of $\Mod R$ and $\mathcal{D}$
be a subcategory of $\Mod R^{op}$ such that $\mathcal{C}^+\subseteq
\mathcal{D}$. Then a homomorphism $f: A\to C$ in $\Mod R$ with $C\in
\mathcal{C}$ is a $\mathcal{C}$-(pre)envelope of $A$ provided $f^+:
C^+\to A^+$ is a $\mathcal{D}$-(pre)cover of $A^+$ in $\Mod
R^{op}$.}

\vspace{0.2cm}

The following are known facts:

(1) $R$ is right coherent and left perfect if and only if every left
$R$-module has a projective preenvelope ([DC, Proposition 3.14] and
[AsM, Proposition 3.5]). A commutative ring $R$ is artinian if and
only if every $R$-module has a projective preenvelope ([AsM,
Corollary 3.6]).

(2) $R$ is right noetherian if and only if every right $R$-module
has an injective (pre)cover ([E, Theorem 2.1]).

So for a right artinian ring $R$, every left $R$-module has a
projective preenvelope and every right $R$-module has an injective
(pre)cover.


We use $\Proj(R)$ (resp. $\Inj(R)$) to denote the subcategory of
$\Mod R$ consisting of projective (resp. injective) left
$R$-modules. As a consequence of Theorem 3.2, we have the following

\vspace{0.2cm}

{\bf Corollary 3.5.} {\it Let $R$ be a left artinian ring and let
$\mathcal{D}$ be a subcategory of $\Mod R^{op}$ containing all
injective modules. Then a homomorphism $f: A\to P$ in $\Mod R$ with
$P$ projective is a projective (pre)envelope of $A$ provided $f^+:
P^+\to A^+$ is a $\mathcal{D}$-(pre)cover of $A^+$ in $\Mod
R^{op}$.}

\vspace{0.2cm}

{\it Proof.} Let $R$ be a left artinian ring. Then every projective
left $R$-module has a decomposition as a direct sum of
indecomposable projective submodules by [AF, Theorem 27.11]. By [SE,
Theorem 5] $\Proj(R)$ is definable, so all projective modules in
$\Mod R$ are pure injective by Lemma 2.2. Note that
$[\Proj(R)]^+\subseteq \Inj(R^{op})$ by Lemma 2.3(1). So the
assertion follows from Theorem 3.2. \hfill{$\square$}

\vspace{0.2cm}

Furthermore we have the following

\vspace{0.2cm}

{\bf Corollary 3.6.} {\it The following statements are equivalent.

(1) $R$ is a left artinian ring.

(2) A monomorphism $f: A\rightarrowtail E$ in $\Mod R$ is an
injective preenvelope of $A$ if and only if $f^+: E^+\to A^+$ is a
projective precover of $A^+$ in $\Mod R^{op}$.

In addition, if $R$ is a left and right artinian ring, then a
homomorphism $f: A\to P$ in $\Mod R$ is a projective preenvelope of
$A$ if and only if $f^+: P^+\to A^+$ is an injective precover of
$A^+$ in $\Mod R^{op}$.}

\vspace{0.2cm}

{\it Proof.} $(2)\Rightarrow (1)$ follows from Lemma 2.3(2).

$(1)\Rightarrow (2)$ It is well known that a right $R$-module is
flat if and only if it is projective over a left artinian ring $R$.
So the assertion follows from [EH, Theorem 3.7].

Let $R$ be a left and right artinian ring. Then
$[\Proj(R)]^+\subseteq \Inj(R^{op})$ and $[\Inj(R^{op})]^+=\Proj(R)$
by Lemma 2.3. Thus the last assertion follows from Corollary 3.5 and
[EH, Corollary 3.2(1)]. \hfill{$\square$}

\vspace{0.2cm}

For a subcategory $\mathcal{C}$ of $\Mod R$, we write
$\mathcal{C}^{\bot}=\{X\in \Mod R\mid \Ext_R^i(C,X)=0$ for any $C\in
\mathcal{C}$ and $i\geq 1\}$ and $^{\bot}\mathcal{C}=\{X\in \Mod
R\mid \Ext_R^i(X,C)=0$ for any $C\in \mathcal{C}$ and $i\geq 1\}$.
We use $\GProj(R)$ (resp. $\GInj(R)$) to denote the subcategory of
$\Mod R$ consisting of Gorenstein projective (resp. injective)
modules. For an artinian algebra $R$, recall from [B1] tha $R$ is
called {\it virtually Gorenstein} if
$[\GProj(R)]^{\bot}={^{\bot}[\GInj(R)}]$, and $R$ is said of {\it
finite Cohen-Macaulay type} ({\it finite CM-type} for short) if
there exist only finitely many non-isomorphic finitely generated
indecomposable Gorenstein projective left $R$-modules. The notion of
virtually Gorenstein algebras is a common generalization of that of
Gorenstein algebras and algebras of finite representation type ([B2,
Example 4.5]).

Note that for a Gorenstein ring (that is, a left and right
noetherian ring with finite left and right self-injective
dimensions) $R$, every finitely generated left $R$-module has a
Gorenstein projective preenvelope ([EJ2, Corollary 11.8.3]).

\vspace{0.2cm}

{\bf Corollary 3.7.} {\it Let $R$ be a virtually Gorenstein artinian
algebra of finite CM-type and let $\mathcal{D}$ be a subcategory of
$\Mod R^{op}$ containing all Gorenstein injective modules. Then a
homomorphism $f: A\to G$ in $\Mod R$ with $G$ Gorenstein projective
is a Gorenstein projective (pre)envelope of $A$ provided $f^+:
G^+\to A^+$ is a $\mathcal{D}$-(pre)cover of $A^+$ in $\Mod
R^{op}$.}

\vspace{0.2cm}

{\it Proof.} Let $R$ be a virtually Gorenstein artinian algebra of
finite CM-type. Then $\GProj(R)$ is definable by [B1], and every
Gorenstein projective module in $\Mod R$ is a direct sum of
indecomposable modules by [B2, Theorem 4.10]. So all Gorenstein
projective modules in $\Mod R$ are pure injective by Lemma 2.2. Note
that $[\GProj(R)]^+\subseteq \GInj(R^{op})$ by [HuX, Corollary 2.6]
and [H, Theorem 3.6]. So the assertion follows from Theorem 3.2.
\hfill{$\square$}

\vspace{0.2cm}

Recall from [Me] that a module $M$ in $\Mod R$ is called {\it
absolutely pure} if it is a pure submodule in every module in $\Mod
R$ that contains it, or equivalently, if it is pure in every
injective module in $\Mod R$ that contains it. Absolutely pure
modules are also known as {\it FP-injective modules}. It is trivial
that an injective module is absolutely pure. By [Me, Theorem 3], a
ring $R$ is left noetherian if and only if every absolutely pure
module in $\Mod R$ is injective.

For a left noetherian ring $R$, every injective left $R$-module has
a decomposition as a direct sum of indecomposable injective
submodules ([M, Theorem 2.5]). It is well known that many results
about finitely generated modules or injective modules over
noetherian rings should have a counterpart about finitely presented
modules or absolutely pure modules (see [G, EJ2, GT, Me, P] and so
on). The following corollary shows that the result just mentioned
above is one of exceptions.

\vspace{0.2cm}

{\bf Corollary 3.8.} {\it For a (left and right coherent) ring $R$,
an absolutely pure left $R$-module does not have a decomposition as
a direct sum of indecomposable absolutely pure submodules in
general.}

\vspace{0.2cm}

{\it Proof.} Let $R$ be a left and right coherent ring. Then the
subcategory of $\Mod R$ consisting of absolutely pure modules is
definable by [K, Proposition 15.1]. If any absolutely pure left
$R$-module has a decomposition as a direct sum of indecomposable
absolutely pure submodules, then any absolutely pure left $R$-module
is pure injective by Lemma 2.2. So by Lemma 2.3(2) and Theorem 3.2,
a homomorphism $f: A\to C$ in $\Mod R$ with $C$ absolutely pure is
an absolutely pure preenvelope of $A$ provided $f^+: C^+\to A^+$ is
a flat precover of $A^+$ in $\Mod R^{op}$, which contradicts [EH,
Example 3.6]. \hfill{$\square$}

\vspace{0.5cm}

{\bf Acknowledgements.} This research was partially supported by the
Specialized Research Fund for the Doctoral Program of Higher
Education (Grant No. 20100091110034), NSFC (Grant No. 11171142), NSF
of Jiangsu Province of China (Grant Nos. BK2010047, BK2010007) and a
Project Funded by the Priority Academic Program Development of
Jiangsu Higher Education Institutions. The author thanks the referee
for the useful suggestions.


\vspace{0.5cm}

\begin{thebibliography}{101}

\bibitem[AF]{A1} F.W. Anderson and K.R. Fuller, Rings and Categories of
Modules, 2nd ed. Grad. Texts in Math. {\bf 13}, Springer-Verlag,
Berlin, 1992.

\bibitem[AsM]{A2} J. Asensio Mayor and J. Mart\'inez Hern\'andez, {\it On flat
projective envelopes}, J. Algebra {\bf 160} (1993), 434--440.

\bibitem[AuS]{A.} M. Auslander and S.O. Smal$\phi$, {\it Preprojective modules over artin algebras},
J. Algebra {\bf 66} (1980), 61--122.

\bibitem[B1]{A4} A. Beligiannis, {\it Cohen-Macaulay modules, (co)torsion pairs,
and virtually Gorenstein algebras}, J. Algebra {\bf 288} (2005),
137--211.

\bibitem[B2]{A5} A. Beligiannis, {\it On algebras of finite Cohen-Macaulay type},
Adv. Math. {\bf 226} (2011), 1973--2019.

\bibitem[CS]{A6} T.J. Cheatham and D. Stone, {\it Flat and projective
character modules}, Proc. Amer. Math. Soc. {\bf 81} (1981),
175--177.

\bibitem[DC]{A7} N.Q. Ding and J.L. Chen, {\it Relative coherence and
preenvelopes}, Manuscripta math. {\bf 81} (1993), 243--262.

\bibitem[E]{A8} E.E. Enochs, {\it Injective and flat
covers, envelopes and resolvents}, Israel J. Math. {\bf 39} (1981),
189--209.

\bibitem[EH]{A9} E.E. Enochs and Z.Y. Huang, {\it Injective envelopes and
(Gorenstein) flat covers}, Algebras Represent. Theory {\bf 15}
(2012), 1131--1145.

\bibitem[EJ1]{A10} E.E. Enochs and O.M.G. Jenda, {\it Gorenstein injective
and projective modules}, Math. Z. {\bf 220} (1995), 611--633.

\bibitem[EJ2]{A11} E.E. Enochs and O.M.G. Jenda, Relative
Homological Algebra, Vol. 1, 2nd Edition, De Gruyter Exp. in Math.
{\bf 30}, Walter de Gruyter, Berlin, Boston, 2011.

\bibitem[F]{A12} D.J. Fieldhouse, {\it Character modules},
Comment. Math. Helv. {\bf 46} (1971), 274--276.

\bibitem[G]{A13} S. Glaz, Commutative Coherent Rings, Lect. Notes in
Math. {\bf 1371}, Springer-Verlag, Berlin, 1989.

\bibitem[GT]{A14} R. G\"obel and J. Trlifaj, Approximations and Endomorphism Algebras of
Modules, de Gruyter Expositions in Math. {\bf 41}, Walter de Gruyter
GmbH \& Co. KG, Berlin, 2006.


\bibitem[H]{A15} H. Holm, {\it Gorenstein homological dimensions}, J. Pure Appl.
Algebra {\bf 189} (2004), 167--193.

\bibitem[HuX]{A16} J.S. Hu and A.M. Xu, {\it On stability of {\bf F}-Gorenstein flat
categories}, Algebra Colloq. (to appear)


\bibitem[K]{A17} H. Krause, {\it The Spectrum of a Module Category}, Memoirs Amer. Math. Soc.
 {\bf 707}, Amer. Math. Soc., Province, RI, 2001.
(2008), 2186--2194.

\bibitem[M]{A18} E. Matlis, {\it Injective modules over noetherian rings}, Pacific J. Math. {\bf 8}
(1958), 511--528.

\bibitem[Me]{A19} C. Megibben, {\it Absolutely Pure Modules}, Proc. Amer. Math. Soc. {\bf 26}
(1970), 561--566.


\bibitem[P]{A20} K. Pinzon, {\it Absolutely pure covers}, Comm. Algebra {\bf 36}
(2008), 2186--2194.

\bibitem[SE]{A21} G. Sabbagh and P. Eklof, {\it Definability problems for rings and
modules}, J. Symbolic Logic {\bf 36} (1971), 623--649.




\end{thebibliography}
\end{document}